\documentclass[10pt]{amsproc}
\usepackage{graphicx,hyperref,color,ytableau}
\newtheorem{theorem}{Theorem}[subsection]

\theoremstyle{definition}
\newtheorem{definition}[theorem]{Definition}
\theoremstyle{remark}

\newtheorem{example}[theorem]{Example}

\newcommand{\rowins}{\mathrm{RINS}}
\newcommand{\ins}{\mathrm{INSERT}}

\title{Knuth's moves on Timed Words}
\author{Amritanshu Prasad}
\date{\today}
\address{The Institute of Mathematical Sciences, Chennai}
\address{Homi Bhabha National Institute, Mumbai}
\begin{document}
\begin{abstract}
  We give an exposition of Schensted's algorithm to find the length of the longest increasing subword of a word in an ordered alphabet, and Greene's generalization of Schensted's results using Knuth equivalence.
  We announce a generalization of these results to timed words.
\end{abstract}
\maketitle
\ytableausetup{smalltableaux}
\section{Introduction}
\label{sec:introduction}
The theory of Young tableaux lies at the cross-roads of modern combinatorics, the theory of symmetric functions, enumerative problems in geometry, and representation theory (see \cite{fulton,manivel,rtcv,schur_poly}).
Young tableaux are named after Alfred Young, who introduced them in his study of the representation theory of symmetric groups, which he called \emph{substitutional analysis} \cite{doi:10.1112/plms/s1-33.1.97}.
Young tableaux played an important role in the proof of the Littlewood-Richardson rule.
This is a rule for computing the Littlewood-Richardson coefficients $c^\lambda_{\mu\nu}$, which arise as (for details, see \cite{fulton,manivel}):
\begin{itemize}
\item the multiplicity of the irreducible polynomial representation $W_\lambda$ of $GL_n(\mathbf C)$ in a tensor product $W_\mu\otimes W_\nu$.
\item the coefficient of the Schur polynomial $s_\lambda$ in the expansion of a product $s_\mu s_\nu$ of Schur polynomials.
\item the number of points  of intersection of Schubert varieties $X_\mu$, $X_\nu$ and $X_{\check\lambda}$ in general position.
\end{itemize}
Robinson \cite{robinson-algo} outlined an approach to the Littlewood-Richardson rule based on Young tableaux, which was perfected in the work of Lascoux and Sch\"utzenberger \cite{plaxique} forty years later.
In this expository article, our point of departure is Schensted's observation \cite{schensted} that Robinson's construction of the  insertion tableau of a word can be used as an algorithm to determine the longest increasing subword of a word in an ordered language.
Schensted used this to give a formula for the number of words with longest increasing subword of a given length.
Schensted's results were generalized by Greene \cite{Greene-schen} using relations which were introduced by Knuth \cite{knuth} to characterize the class of words with a given insertion tableau.

In a different context, Alur and Dill \cite{alur-dill} introduced timed words as a part of their description of timed automata.
Timed automata are generalizations of finite automata, and are used for the formal verification of real-time systems.
The author has extended Greene's theorem to timed words \cite{timed-plactic}, with the goal of providing a framework to study piecewise linear versions of bijective correspondences involving Young tableau, such as the ones studied by Berenstein and Kirillov \cite{kir-trop}.
The salient features of this extension are outlined here.
Detailed proofs, technical details, and applications to piecewise-linear bijections will appear in \cite{timed-plactic}.
\section{Schensted's Algorithm}
\label{sec:intro}
\subsection{Words}
\label{sec:words}
Let $A_n$ denote the set $\{1,\dotsc,n\}$, which we regard as an \emph{ordered alphabet}.
A \emph{word} in $A_n$ is a finite sequence $w=c_1\dotsb c_k$ of elements of $A_n$.
The set of all words in $A_n$ is denoted $A_n^*$.
A \emph{subword} of $w=c_1\dotsb c_k$ is defined to be a word of the form
\begin{displaymath}
  w' = c_{i_1}\dotsb c_{i_m}, \text{ where } 1\leq i_1 < \dotsb < i_m \leq k.
\end{displaymath}
The subword $w'$ is said to be \emph{weakly increasing} if $c_{i_1}\leq \dotsb \leq c_{i_m}$.

Consider the following computational problem:
\begin{quote}
  Given a word $w\in A_n^*$, determine the maximal length of a weakly increasing subword of $w$.
\end{quote}
\subsection{Tableaux}
Schensted~\cite{schensted} gave an elegant algorithm to solve the preceding computational problem.
His algorithm makes one pass over the word.
At each stage of its running, it stores a combinatorial object called a \emph{semistandard Young tableau} (see Section~\ref{definition:ssyt}).
This tableau is modified as each successive letter of the word is read.
The length of the longest increasing subword can be read off from the tableau (see Sections~\ref{sec:insert-tabl-word} and~\ref{sec:schensted-theorem}) obtained when all of $w$ has been read.

\begin{definition}
  [Semistandard Young Tableau]
  \label{definition:ssyt}
  A semistandard Young \linebreak tableau in $A_n$ is a finite arrangement of integers from $A_n$ in rows and columns so that the numbers increase weakly along rows, strictly along columns, so that there is an element in the first row of each column, there is an element in the first column of each row, and there are no gaps between numbers.

  Let $l$ be the number of rows in the tableau, and for each $i=1,\dotsc,l$, let $\lambda_i$ be the length of the $i$th row.
  Then $\lambda=(\lambda_1,\dotsc,\lambda_l)$ is called the \emph{shape} of the tableau.
\end{definition}
\begin{example}
  \label{example:ssyt}
  The arrangement
  \begin{displaymath}
    \ytableaushort{115,24,3}
  \end{displaymath}
  is a semistandard Young tableau of shape $(3,2,1)$ in $A_5$.
\end{example}

The notion of a semistandard Young tableau is a generalization of \emph{Young tableau}, which was introduced by Young~\cite[p.~133]{doi:10.1112/plms/s1-33.1.97}.
In Young's version, each element of $A_n$ occurs exactly once in the tableau.
For brevity, we shall henceforth use the term \emph{tableau} to refer to a semistandard Young tableau.
\subsection{Row Insertion}
\label{sec:row-insertion}
A word $c_1c_2\dotsb c_k$ in $A_n^*$ is called a \emph{row} if $c_1\leq \dotsb \leq c_k$.
Each row of a tableau is a row in the sense of this definition.
For each row $u=a_1\dotsb a_k\in A_n^*$, define the \emph{row insertion} of $a$ into $u$ by:
\begin{displaymath}
  \rowins(u,a) =
  \begin{cases}
    (\emptyset, a_1\dotsb a_k a) & \text{if } a_k\leq a,\\
    (a_j,a_1\dotsb a_{j-1}aa_{j+1}\dotsb a_k) & \text{otherwise, with}\\
    & j=\min\{i\mid a<a_i\}.
  \end{cases}
\end{displaymath}
Here $\emptyset$ should be thought of as an empty word of length zero.
\begin{example}
  $\rowins(115,5) = (\emptyset,1155)$, $\rowins(115,3)=(5,113)$.
\end{example}
It is clear from the construction that, for any row $u\in A_n^*$ and $a\in A_n$, if $(a',u')=\rowins(u,a)$, then $u'$ is again a row.
For convenience set $\rowins(u,\emptyset)=(\emptyset,u)$.
\subsection{Tableau Insertion}
\label{sec:tableau-insertion}
Let $t$ be a tableau with rows $u_1,u_2,\dotsc, u_l$.
Then $\ins(t,a)$, the insertion of $a$ into $t$, is defined as follows: first $a$ is inserted into $u_1$; if $\rowins(u_1,a)=(a_1',u_1')$, then $u_1$ is replaced by $u_1'$.
Then $a_1'$ is inserted into $u_2$; if $\rowins(u_2,a_1')=(a_2',u_3)$, then $u_2$ is replaced by $u_2'$, $a_2'$ is inserted into $u_3$, and so on.
This process continues, generating $a_1',a_2',\dotsc,a_k'$ and $u_1',\dotsc,u_k'$.
The tableau $t'=\ins(t,a)$ has rows $u_1',\dotsc,u_k'$, and a last row (possibly empty) consisting of $a_k'$.
It turns out that $\ins(t,a)$ is a tableau \cite{knuth}.
\begin{example}
  \label{example:insertion}
  For $t$ as in Example~\ref{example:ssyt}, we have
  \begin{displaymath}
    \ins(t,3) = \ytableaushort{113,245,3},
  \end{displaymath}
  since $\rowins(115,3)=(5,113)$, $\rowins(24,5)=(\emptyset,245)$.
\end{example}
\subsection{Insertion Tableau of a Word}
\label{sec:insert-tabl-word}
\begin{definition}
\label{definition:insertion-tableau}
The insertion tableau $P(w)$ of a word $w$ is defined recursively as follows:
\begin{align}
  P(\emptyset)&=\emptyset\\
  P(c_1\dotsb c_k)&=\ins(P(c_1\dotsb c_{k-1}), c_k).
\end{align}
\end{definition}
\begin{example}
  \label{example:insertion-tableau}
  Take $w=3421153$.
  Sequentially inserting the terms of $w$ into the empty tableau $\emptyset$ gives the sequence of tableaux:
  \begin{displaymath}
    \ytableaushort{3},\ytableaushort{34},\ytableaushort{24,3},\ytableaushort{14,2,3},\ytableaushort{11,24,3},\ytableaushort{115,24,3},
  \end{displaymath}
  and finally, the insertion tableau $P(w)=\ytableaushort{113,245,3}$.
\end{example}
\subsection{Schensted's Theorem}
\label{sec:schensted-theorem}
Schensted \cite{schensted} proved the following:
\begin{theorem}
  The length of the longest increasing subword of any $w\in A_n^*$ is the length of the first row of $P(w)$.
\end{theorem}
In other words, the algorithm for constructing the insertion tableau of $w$ solves the computational problem posed in Section~\ref{sec:words}.

The proof of Schensted's theorem is not very difficult, and the reader is invited to attempt it.
The proof is by induction on $k$, and uses the observation is that the last entry of the first row of $P(a_1\dotsb a_k)$ is the \emph{least last element} of all maximal length weakly increasing subword of $a_1\dotsb a_k$.
\subsection{Greene's Theorem}
\label{sec:greenes-theorem}
The insertion tableau $P(w)$ obtained from a word $w$ seems to contain a lot more information than just the length of the longest weakly increasing subword.
For example, what do the lengths of the remaining rows of $P(w)$ signify?
The answer to this question was given by Greene~\cite{Greene-schen}.
We say that subwords $c_{1_1}\dotsb c_{i_r}$ and $c_{j_1}\dotsb c_{j_s}$ of $c_1\dotsb c_k$ are \emph{disjoint} if the subsets $\{i_1,\dotsc,i_r\}$ and $\{j_1,\dotsc,j_s\}$ are disjoint.
\begin{definition}
  [Greene Invariants]
  The $r$th Greene invariant of a word $w\in A_n^*$ is defined to be the maximum cardinality of a union of $r$ pairwise disjoint weakly increasing subwords of $w$.
\end{definition}
\begin{example}
  \label{example:greene}
  For $w=3421153$ from Example~\ref{example:insertion-tableau}, the longest weakly increasing subwords have length $3$ (for example, $113$ and $345$).
  The subwords $345$ and $113$ are disjoint, and no pair of disjoint weakly increasing subwords of $w$ can have cardinality greater than $6$.
  However, the entire word $w$ is a union of three disjoint weakly increasing subwords (for example $345$, $23$ and $15$).
  So the Greene invariants of $w$ are $a_1(w)=3$, $a_2(w)=6$, and $a_3(w)=7$.
\end{example}
\begin{theorem}
  [Greene]
  \label{theorem:greene}
  For any $w\in A_n$, if $P(w)$ has shape $\lambda=(\lambda_1,\dotsc, \lambda_l)$, then for each $r=1,\dotsc,l$, $a_r(w)=\lambda_1+\dotsb + \lambda_l$.
\end{theorem}
Example~\ref{example:greene} is consistent with Greene's theorem as the shape of $P(w)$ is $(3, 3, 1)$ and the Greene invariants are $3$, $6=3+3$ and $7=3+3+1$, respectively.

\subsection{Knuth Equivalence}
\label{sec:knuth-equivalence}

Greene's proof of Theorem~\ref{theorem:greene} is based on the notion of Knuth equivalence.
Knuth \cite{knuth} identified a pair of elementary moves on words:
\begin{gather}
  \tag{$K1$}\label{eq:k1}
  xzy \equiv zxy \text{ if } x\leq y < z,
  \\
  \tag{$K2$}\label{eq:k2}
  yxz \equiv yzx \text{ if } x < y \leq z.
\end{gather}
For example, in the word $4213443$, the segment $213$ is of the form $yxz$, with $x<y\leq z$.
A Knuth move of type (\ref{eq:k2}) replaces this segment by $yzx$, which is $231$.
Thus a Knuth move of type (\ref{eq:k2}) transforms $4\textbf{213}443$ into $4\textbf{231}443$.
Knuth equivalence is the equivalence relation on $A_n^*$ generated by Knuth moves:
\begin{definition}[Knuth Equvalence]
  \label{definition:Knuth-equiv}
  Words $w,w'\in A_n^*$ are said to be Knuth equivalent if $w$ can be transformed into $w'$ by a series of Knuth moves (\ref{eq:k1}) and (\ref{eq:k2}).
  If this happens, we write $w\equiv w'$.
\end{definition}
\begin{example}
  \label{example:knuth-red}
  The word $3421153$ is Knuth equivalent to $3245113$:
  \begin{displaymath}
    \textbf{342}1153 \equiv_{K2} 3\textbf{241}153 \equiv_{K2} 32\textbf{141}53 \equiv_{K1} 3241\textbf{153} \equiv_{K1} 324\textbf{151}3 \equiv_{K1} 3245113.
  \end{displaymath}
  At each stage, the letters to which the Knuth moves will be applied to obtain the next stage are highlighted.
\end{example}
\subsection{Reading Word of a Tableau}
\label{sec:reading-word}
Given a tableau, its reading word is obtained by reading its rows from left to right, starting with the bottom row, and moving up to its first row.
\begin{example}
  \label{example:reading-word}
  The reading word of the tableau:
  \begin{displaymath}
    \ytableaushort{113,245,3}
  \end{displaymath}
  is $3245113$.
\end{example}
\subsection{Proof of Greene's Theorem}
\label{sec:proof-greene}
The proof of Greene's theorem is based on three observations, all fairly easy to prove:
\begin{enumerate}
\item If $w$ is the reading word of a tableau of shape $\lambda=(\lambda_1,\dotsc,\lambda_l)$, then $a_r(w)=\lambda_1+\dotsb + \lambda_r$ for $r=1,\dotsc,l$.
\item Every word is Knuth equivalent to the reading word of its insertion tableau.
\item Greene invariants remain unchanged under Knuth moves.
\end{enumerate}
We illustrate these points with examples (for detailed proofs, see Lascoux, Leclerc and Thibon~\cite{Lascoux}, or Fulton \cite{fulton}).
For the first point, in Example~\ref{example:reading-word} the first $k$ rows of the tableau
\begin{displaymath}
  \ytableaushort{113,245,3}
\end{displaymath}
are indeed disjoint weakly increasing subwords of its reading word of maximal cardinality.
For the second point, observe that the sequence of Knuth moves in Example~\ref{example:knuth-red} transform $3421153$ to the reading word of its insertion tableau.
For the third point, consider the case of the Knuth move (\ref{eq:k1}).
A word of the from $w=uxzyv$ is transformed into the word $w'=uzxyv$.
The only issue is that a weakly increasing subword $g$ of $w$ may contain both the letters $x$ and $z$.
Then it no longer remains a weakly increasing subword of $w'$.
However, the subword, being weakly increasing, cannot contain $y$, so the $z$ can be swapped for a $y$.
This could be a problem if $y$ is part of another weakly increasing subword $g'$ in a collection of pairwise disjoint weakly increasing subwords.
In that case, we have $g=g_1xzg_2$ and $g'=g'_1y g'_2$.
We may replace them with $g_1xyg'_2$ and $g_1'zg_2$, which would still be weakly increasing, and would have the same total length as $g$ and $g'$.
\subsection{Characterization of Knuth Equivalence}
\label{sec:characterization}
Knuth equivalence can be characterized in terms of Greene invariants (see \cite[Theorem~2.15]{plaxique}).
\begin{theorem}
  Two words $w$ and $w'$ in $A_n^*$ are Knuth equivalent if and only if $a_r(uwv)=a_r(uw'v)$ for all words $u$ and $v$ in $A_n^*$, and all $r\geq 1$.
\end{theorem}
\section{Timed Words}
\subsection{From Words to Timed Words}
\label{sec:timed-words}
\label{sec:words-to-timed-words}
Words, in the sense of Section~\ref{sec:intro}, play an important role in computer science, specifically in the formal verification of systems.
Each letter of the alphabet is thought of as an event.
A sequence of events it then nothing but a word in $A_n^*$.
The system is modeled as an automaton having a starting state, and each time an event occurs, its state changes, depending both, on its current state, and the event that has occurred.
Following the groundbreaking work of Rabin and Scott \cite{rabin1959finite}, finite state automata are widely used to model and formally verify the integrity of systems.

For many real-time systems, such as controllers of washing machines, industrial processes, and air or railway traffic control, the \emph{time gaps} between the occurrences of the events modeled by words are as important as the events themselves.

To deal with real-time systems, Alur and Dill \cite{alur-dill} developed the theory of \emph{timed automata}.
A timed automaton responds to a sequence of events that come with time stamps for their occurrence.
They represented a sequence of events with time stamps by timed words.
We introduce a finite variant of the notion of timed word that they used:
\begin{definition}
  [Timed Word]
  \label{definition:timed-word}
  A timed word in $A_n$ is a sequence of the form:
  \begin{equation}
    \label{eq:timed-word}
    w=c_1^{t_1} c_2^{t_2}\dotsb c_k^{t_k},
  \end{equation}
  where $c_1,\dotsc,c_k\in A_n$, and $t_1,\dotsc,t_k$ are positive real numbers, and $c_i\neq c_{i+1}$ for $i=1,\dotsc,k-1$.
  The length of the timed word $w$ above is $l(w)=t_1+\dotsb+t_k$.
\end{definition}
A sequence (\ref{eq:timed-word}) where $c_i=c_{i+1}$ also represents a timed word; segments of the form $c^{t_1}c^{t_2}$ are replaced by $c^{t_1+t_2}$ until all consecutive pairs of terms have different letters.
The timed word $w$ in (\ref{eq:timed-word}) may also be regarded as a piecewise constant left-continuous function $\mathbf w:[0,l(w))\to A_n$, where
\begin{displaymath}
  \mathbf w(t) = c_i \text{ if } t_1+\dotsb+t_{i-1}\leq t < t_1+\dotsb + t_i.
\end{displaymath}
The function $\mathbf w:[0, l(w))\to A_n$ is called the \emph{function associated to the timed word $w$}.
We say that the timed word $w$ is a \emph{timed row} if $c_1<\dotsb < c_k$.
Timed words form a monoid under concatenation.
The monoid of timed words in the alphabet $\{1,\dotsc,n\}$ is denoted $A_n^\dagger$.
The map:
\begin{displaymath}
  a_1\dotsb a_k \mapsto a_1^1 a_2^1 \dotsb a_k^1
\end{displaymath}
defines an embedding of $A_n^*$ in $A_n^\dagger$ as a submonoid.
\begin{example}
  \label{example:timed-word}
  An example of a timed word in $A_6^\dagger$ of length $7.19$ is:
  \begin{displaymath}
    w = 3^{0.82}5^{0.08}2^{0.45}6^{0.64}5^{0.94}1^{0.15}5^{0.09}1^{0.52}4^{0.29}1^{0.59}3^{0.97}4^{0.42}2^{0.61}1^{0.07}4^{0.55}
  \end{displaymath}
  Using a color-map to represent the integers $1$ to $6$,
  \begin{center}
    \includegraphics[width=0.85\textwidth]{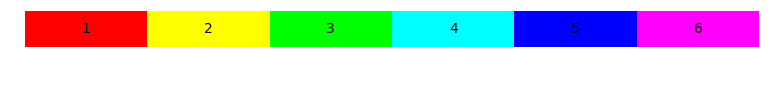}
  \end{center}
  the timed word $w$ can be visualized as a colored ribbon:
  \begin{center}
    \includegraphics[width=\textwidth]{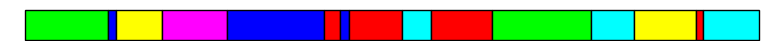}.
  \end{center}
\end{example}
\subsection{Subwords of Timed Words}
\label{sec:subwords-timed-words}
\begin{definition}
  [Time Sample]
  A \emph{time sample} of a word $w$ is a subset of $[0,l(w))$ of the form:
  \begin{displaymath}
    S = [a_1,b_1)\cup \dotsb \cup [a_k,b_k),
  \end{displaymath}
  where $0\leq a_1 < b_1 < a_2 < b_2 < \dotsb <a_k < b_k \leq l(w)$.
  The length of the time sample $S$ is $\sum_i (b_i-a_i)$, the Lebesgue measure $\mu(S)$ of $S$.
\end{definition}
Given a time sample $S\subset [0, l(w))$, and $0\leq t\leq l(S))$, the set
\begin{displaymath}
  \{\tilde t \mid \mu(S\cap [0,\tilde t)) = t\}
\end{displaymath}
is a closed interval $[a_t,b_t]\subset [0, l(S))$.
This happens because the function $t'\mapsto \mu(S\cap [0,t'))$ is a piecewise-linear continuous function on $[0,l(w)]$ which takes value $0$ at $t'=0$, and $l(S)$ at $t'=1$.
\begin{definition}
  [Subword of a Timed Word]
  \label{definition:timed-subword}
  The subword of a timed word with respect to a time sample $S\subset [0, l(w))$ is the timed word $w_S$ of length $\mu(S)$ whose associated function is given by:
  \begin{displaymath}
    \mathbf w_S(t) =  \mathbf w(b_t) \text{ for } 0\leq t < \mu(S),
  \end{displaymath}
  where $b_t$ is the largest number in $[0,l(w))$ such that $\mu(S\cap [0, \tilde t)) = t$.
\end{definition}
\subsection{Timed Tableau}
\label{sec:Timed-Tableau}
\begin{definition}
  [Timed Tableau]
  \label{definition:timed-tableau}
  A timed tableau is a collection $u_1,u_2,\dotsc, u_l$ of timed words such that
  \begin{enumerate}
  \item Each $u_i$ is a timed row (in the sense of Section~\ref{sec:timed-words}).
  \item For each $i=1,\dotsc,l-1$, $l(u_i)\geq l(u_{i+1})$.
  \item For each $i=1,\dotsc,l-1$ and $0\leq t<l(u_{i+1})$, $u_i(t)<u_{i+1}(t)$.
  \end{enumerate}
\end{definition}
\begin{example}
  \label{example:timed-tableau}
  A timed tableau of shape $(3.20,1.93,1.09,0.61,0.29,0.07)$ is given by:
  \begin{align*}
    t = & 1^{1.33}2^{0.54}3^{0.36}4^{0.97}\\
    & 2^{0.52}3^{0.91}5^{0.50}\\
    &3^{0.52}4^{0.22}5^{0.32}6^{0.03}\\
    &4^{0.07}5^{0.22}6^{0.32}\\
    &5^{0.07}6^{0.22}\\
    &6^{0.07}
  \end{align*}
  In using the color-map from Section~\ref{sec:timed-words}, it can be visualized as:
  \begin{center}
    \includegraphics[width=0.5\textwidth]{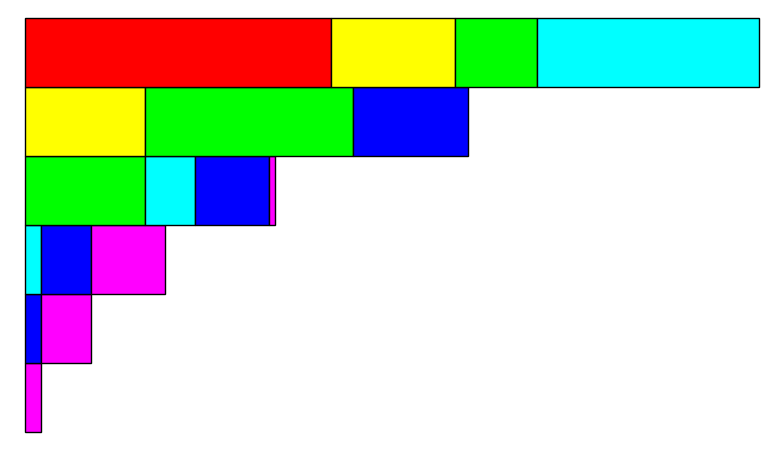}
  \end{center}
  The three properties of Definition \ref{definition:timed-tableau} are easily perceived from the figure.
\end{example}
\begin{definition}
  [Reading Word of a Timed Tableau]
  The reading word of a timed tableau with rows $u_1,\dotsc,u_l$ is the timed word
  \begin{displaymath}
    u_l u_{l-1}\dotsb u_1.
  \end{displaymath}
\end{definition}
\begin{example}
  The reading word of the timed tableau in Example~\ref{example:timed-tableau} is
  \begin{displaymath}
    6^{0.07}\;5^{0.07}6^{0.22}\;4^{0.07}5^{0.22}6^{0.32}\;3^{0.52}4^{0.22}5^{0.32}6^{0.03}\;2^{0.52}3^{0.91}5^{0.50}\;1^{1.33}2^{0.54}3^{0.36}4^{0.97}.
  \end{displaymath}
\end{example}
\subsection{Timed Insertion}
\label{sec:timed-insertion}
Given a timed word $w$ and $0\leq a < b \leq l(w)$, according to Definition~\ref{definition:timed-subword}, $w_{[a, b)}$ is the timed word of length $b-a$ such that:
\begin{displaymath}
  w_{[a, b)}(t) = w(a+ t) \text{ for } 0\leq t<b-a.
\end{displaymath}
\begin{definition}[Timed Row Insertion]
  \label{definition:timed-row-insertion}
  Given a timed row $w$, define the insertion $\rowins(w, c^{t_c})$ of $c^{t_c}$ into $w$ as follows: if $w(t)\leq c$ for all $0\leq t < l(u)$, then
  \begin{displaymath}
    \rowins(w, c^{t_c}) = (\emptyset, wc^{t_c}).
  \end{displaymath}
  Otherwise, there exists $0\leq t < l(w)$ such that $\mathbf w(t)>c$.
  Let
  \begin{displaymath}
    t_0 = \min\{0\leq t< l(w)\mid \mathbf w(t)> c\}.
  \end{displaymath}
  Define
  \begin{displaymath}
    \rowins(w, c^{t_c}) =
    \begin{cases}
      (w_{[t_0, t_0+t_c)}, w_{[0, t_0)}c^{t_c} w_{[t_0+t_c, l(w))}) & \text{if } l(w) - t_0 > t_c,\\
      (w_{[t_0, l(u))}, w_{[0, t_0)} c^{t_c}) & \text{if } l(w) - t_0 \leq t_c.
    \end{cases}
  \end{displaymath}
  It is obvious that the above definition is compatible with the definition of $\rowins$ from Section~\ref{sec:tableau-insertion} when $u$ is a row in $A_n^*$, and $t_c=1$.
  If $u=c_1^{t_1}\dotsb c_l^{t_l}$ is a timed word, define $\rowins(w,u)$ by induction on $l$ as follows:
  Having defined $(v',w')=\rowins(w,c_1^{t_1}\dotsb c_{l-1}^{t_{l-1}})$,
  let $(v'',w'')=\rowins(w',c_l^{t_l})$.
  Then define
  \begin{displaymath}
    \rowins(w,u) = (v'v'', w'').
  \end{displaymath}
\end{definition}
\begin{example}
  \label{example:timed-row-ins}
  $\rowins(1^{1.4}2^{1.6}3^{0.7},1^{0.7}2^{0.2})=(2^{0.7}3^{0.2},1^{2.1}2^{1.1}3^{0.5})$.
\end{example}
\begin{definition}
  [Timed Tableau Insertion]
  Let $w$ be a timed tableau with row decomposition $u_l\dotsc u_1$, and let $v$ be a timed row.
  Then $\ins(w, v)$, the insertion of $v$ into $w$ is defined as follows:
  first $v$ is inserted into $u_1$.
  If $\rowins(u_1,v)=(v_1',u_1')$, then $v_1'$ is inserted into $u_2$; if $\rowins(u_2,v_1')=(v_2',u_2')$, then $v_2'$ is inserted in $u_3$, and so on.
  This process continues, generating $v_1',\dotsc,v_l'$ and $u_1',\dotsc,u_l'$.
  $\ins(t,v)$ is defined to be $v_l'u_l'\dotsb u_1'$.
  Note that it is quite possible that $v_l'=\emptyset$.
\end{definition}
\begin{example}
  Take
  \begin{align*}
    w = & 1^{1.4}2^{1.6}3^{0.7}\\
    & 3^{0.8}4^{1.1},
  \end{align*}
  a timed tableau in $A_5$ of shape $(3.7,1.9)$.
  Then
  \begin{align*}
    \ins(w,1^{0.7}2^{0.2})=& 1^{1.7}2^{3}3^{0.2}\\
    & 2^{0.3}3^{1.2}4^{0.4}\\
    & 3^{0.3}4^{0.7}
  \end{align*}
  of shape $(4.9,1.9,1.0)$.
\end{example}
\subsection{Insertion Tableau of a Timed Word}
\label{sec:insertion-tableau}
\begin{definition}
  [Insertion Tableau of a Timed Word]
  The insertion tableau $P(w)$ of a timed word $w$ is defined recursively by the rules:
  \begin{enumerate}
  \item $P(\emptyset) = \emptyset$,
  \item $P(wc^t) = \ins(P(w),c^t)$.
  \end{enumerate}
\end{definition}
\begin{example}
  The tableau in Example~\ref{example:timed-tableau} is the insertion tableau of the timed word in Example~\ref{example:timed-word}.
\end{example}
\subsection{Greene Invariants for Timed Words}
\label{sec:greene-invar-timed}
\begin{definition}
  [Greene Invariants for Timed Words]
  The $r$th Greene invariant for a timed word $w$ is defined as:
  \begin{displaymath}
    a_r(w) = \sup\left\{\mu(S_1)+\dotsb + \mu(S_r)\left|
      \begin{array}{cc}
        S_1,\dotsc, S_r \text{ are pairwise disjoint time samples}\\
        \text{of $w$ such that $w_{S_i}$ a timed row for each $i$}
      \end{array}
      \right.
    \right\}.
  \end{displaymath}
\end{definition}
\subsection{Greene's Theorem for Timed Words}
\label{sec:greene-theorem-timed}
All the ingredients are now in place to state Greene's theorem for timed words:
\begin{theorem}
  [Greene's Theorem for Timed Words]
  \label{theorem:Greeene-timed}
  Let $w\in A_n^\dagger$ be a timed word.
  Suppose that $P(w)$ has shape $\lambda=(\lambda_1,\dotsc,\lambda_l)$, then the Greene invariants of $w$ are given by:
  \begin{displaymath}
    a_r(w) = \lambda_1+\dotsb + \lambda_r \text{ for } r=1,\dotsc,l.
  \end{displaymath}
  For the word $w$ from Example~\ref{example:timed-word}, the insertion tableau has shape
  \begin{displaymath}
    (3.20, 1.93, 1.09, 0.61, 0.29, 0.07),
  \end{displaymath}
  (given in Example~\ref{example:timed-tableau}) so the Greene invariants are given by:
  \begin{align*}
    a_1(w) & = 3.20\\
    a_2(w) & = 3.20+1.93=5.13\\
    a_3(w) & = 3.20+1.93+1.09=6.22\\
    a_4(w) & = 3.20+1.93+1.09+0.61=6.83\\
    a_5(w) & = 3.20+1.93+1.09+0.69+0.29=7.12\\
    a_6(w) & = 3.20+1.93+1.09+0.69+0.29+0.07=7.19
  \end{align*}
\end{theorem}
\subsection{Knuth Moves on Timed Words}
\label{sec:knuth-moves-timed}
As explained in Section~\ref{sec:proof-greene}, the proof of Greene's theorem in \cite{Greene-schen} uses Knuth moves to reduce to the case of reading words of tableau.
The main difficulty in generalizing his theorem to timed words is to identify the analogues of Knuth relations (\ref{eq:k1}) and (\ref{eq:k2}).
These relations need to be simple enough so that it can we shown that if two words differ by such a relation, then they have the same Knuth invariants.
At the same time, they need to be strong enough to reduce any timed word to its insertion tableau.

Consider the relations:
\begin{align*}
  \tag{$\kappa_1$}
  \label{eq:tk1}
  xzy & \equiv zxy \text{ when $xyz$ is a timed row, $l(z)=l(y)$, and $\lim_{t\to l(y)^-}\mathbf y(t)<\mathbf z(0)$},\\
  \tag{$\kappa_2$}
  \label{eq:tk2}
  yxz & \equiv yzx \text{ when $xyz$ is a timed row, $l(x)=l(y)$, and $\lim_{t\to l(x)^-}\mathbf x(t)<\mathbf y(0)$}.
\end{align*}
\begin{example}
  We have:
  \begin{displaymath}
    w = 5^{1.10}3^{2.19}4^{0.89}5^{1.20}1^{0.32}2^{0.44}\equiv  w' = 5^{1.10}3^{2.19}4^{0.62}1^{0.32}2^{0.41}4^{0.27}5^{1.20}2^{0.03},
  \end{displaymath}
  because we may write
  \begin{align*}
    w & = 5^{1.10}3^{2.08}yzx2^{0.03},\\
    w' & = 5^{1.10}3^{2.08}yxz2^{0.03},
  \end{align*}
  where $x=1^{0.32}2^{0.41}$, $y=3^{0.11}4^{0.62}$, and $z=4^{0.27}5^{1.20}$, so $w$ and $w'$ differ by a Knuth move of the form (\ref{eq:tk2}).
\end{example}
We say that two timed words $w$ and $w'$ are Knuth equivalent (denoted $w\equiv w'$) if $w$ can be obtained from $w'$ by a sequence of Knuth moves of the form (\ref{eq:tk1}) and (\ref{eq:tk2}).

With these definitions, we have the following results, which suffice to complete the proof of Theorem~\ref{theorem:Greeene-timed}:
\begin{enumerate}
\item if $w$ is the reading word of a timed tableau of shape $\lambda=(\lambda_1,\dotsc,\lambda_l)$, then $a_r(w)=\lambda_1+\dotsb+\lambda_r$ for $r=1,\dotsc,l$.
\item for every $w\in A_n^\dagger$, $w\equiv P(w)$,
\item if $w\equiv w'$, then $a_r(w)=a_r(w')$ for all $r$.
\end{enumerate}
\subsection{Characterization of Knuth Equivalence}
\label{sec:char-knuth-equiv-timed}
Finally, it turns out that Knuth equivalence for timed words is characterized by Greene invariants, just as in the classical setting (Section~\ref{sec:characterization}):
\begin{theorem}
  Given timed words $w, w'\in A_n^\dagger$, $w\equiv w'$ if and only if, for all $u, v\in A_n^\dagger$,
  \begin{displaymath}
    a_r(uwv) = a_r(uw'v) \text{ for all } r>0.
  \end{displaymath}
\end{theorem}

\end{document}